 \documentclass[preprint,review,12pt]{elsarticle}
\usepackage{graphicx}

\usepackage{amssymb,amsmath}
\usepackage{amsthm}

\journal{}

\begin{document}

\begin{frontmatter}



\title{Effects of the noise level on stochastic semilinear fractional heat equations}

\author[]{Kexue Li}
\ead{kexueli@gmail.com}

\address{School of Mathematics and Statistics, Xi'an Jiaotong University, Xi'an 710049, China}

\begin{abstract}
We consider the stochastic fractional heat equation $\partial_{t}u=\triangle^{\alpha/2}u+\lambda\sigma(u)\dot{w}$ on $[0,L]$ with Dirichlet boundary conditions, where $\dot{w}$  denotes the space-time white noise. For any $\lambda>0$, we prove that the $p$th moment of $\sup_{x\in [0,L]}|u(t,x)|$ grows at most exponentially. Moreover, we prove that the $p$th moment of $\sup_{x\in [0,L]}|u(t,x)|$ is exponentially stable if $\lambda$ is small. At last, We obtain the noise excitation index of $p$th energy of $u(t,x)$ as $\lambda\rightarrow \infty$.

\end{abstract}


\end{frontmatter}
\section{Introduction}
We consider the following stochastic semilinear fractional heat equation with Dirichlet boundary condition:
\begin{equation}\label{sfh}
 \ \left\{\begin{aligned}
&\partial_{t}u(t,x)=\triangle^{\alpha/2}u(t,x)+\lambda\sigma(u(t,x))\dot{w}(t,x),\ t>0,\ x\in (0,L),\\
& u(t,x)=0,\ t>0, \ x\in (0,L)^{c}\\
&u(0,x)=u_{0}(x), \ x\in (0,L),
\end{aligned}\right.
\end{equation}
where $L$ is a positive constant, $(0,L)^{c}=\mathbb{R}\backslash (0,L)$,  $\lambda$ is a positive parameter, $\dot{w}$ is the space-time white noise on $(0,\infty)\times [0,L]$, $\sigma: \mathbb{R}\rightarrow \mathbb{R}$ is a continuous function, $\triangle^{\alpha/2}:=-(-\Delta)^{\alpha/2}$ denotes the fractional Laplacian defined by
\begin{equation}
\triangle^{\alpha/2}f(x)=c(\alpha)\lim_{\varepsilon\downarrow 0}\int_{|y|>\varepsilon}\frac{f(x+y)-f(x)}{|y|^{1+\alpha}}dy,
\end{equation}
where $y\in \mathbb{R}$, $1<\alpha< 2$, $c(\alpha)$ is a positive constant that depends only on $\alpha$. It is known that $\Delta^{\alpha/2}$ is the $L^{2}$ generator
of a rotationally symmetric $\alpha$-stable process, see \cite{CMN}. \\
The initial value $u_{0}$ satisfies the following condition. \\
\textbf{Assumption 1.1.} $u_{0}$ is non-random and continuous on $[0,L]$, $\inf_{x\in [\mu,L-\mu]}u_{0}(x)>0$, where $supp(u_{0})$ denotes the support of $u_{0}$ and $\mu\in (0,L/2)$.\\
The function $\sigma$ satisfies the following condition. \\
\textbf{Assumption 1.2.} There exist constants $l_{\sigma}$ and $L_{\sigma}$ such that for $u,v\in \mathbb{R}$,
$|\sigma(u)-\sigma(v)|\leq L_{\sigma}|u-v|$ and $l_{\sigma}|u|\leq |\sigma(u)|\leq L_{\sigma}|u|$.

We denote by $\{\mathcal{F}_{t}\}_{t\geq 0}$ the filtration generated by the Brownian sheet $\{w(t,x); t\geq 0, x\in [0,L]\}$.
For $1\leq p<\infty$, $L^{p}[0,L]$ denotes the space of all real valued measurable functions $u: [0,L]\rightarrow \mathbb{R}$, such that $|u(x)|^{p}$ is integrable for $x\in [0,L]$. It is a Banach space when equipped with by $\|u\|_{L^{p}[0,L]}=\big(\int_{0}^{L}|u(x)|^{p}dx\big)^{1/p}$. If $p=\infty$, the space
$L^{p}[0,L]$ consists of all real valued measurable functions with a finite norm $\|u\|_{L^{\infty}}=\mbox{ess} \sup_{x\in[0,L]}|u(x)|$.
For each $\gamma\in \mathbb{R}$ and $p\geq 2$, by $B_{p,\gamma}$ we denote the class of all the $\mathcal{F}_{t}$-adapted and continuous random fields $\{u(t,x), t\geq 0, x\in [0,L]\}$ such that $\sup_{t\geq 0}\mathbb{E}[e^{\gamma t}\|u(t)\|_{L^{\infty}}^{p}]<\infty$. It is a Banach space when normed by
\begin{align}\label{gamma}
\|u\|_{p,\gamma}=\big(\sup_{t\geq 0}\mathbb{E}[e^{\gamma t}\|u(t)\|_{L^{\infty}}^{p}]\big)^{1/p}.
\end{align}

Following the approach of Walsh \cite{Walsh}, the $\mathcal{F}_{t}$-adapted random field $\{u_{\lambda}(t,x), t\geq 0, x\in [0,L]\}$ is a mild solution of \emph{Eq.} (\ref{sfh}) with initial value $u_{0}$ and the Dirichlet boundary condition if the following integral equation holds
\begin{align}\label{Walsh}
u_{\lambda}(t,x)=\int_{0}^{L}p_{D}(t,x,y)u_{0}(y)dy+\int_{0}^{t}\int_{0}^{L}p_{D}(t-s,x,y)\lambda\sigma(u(s,y))w(dsdy)
\end{align}
where $p_{D}(t,x,y)$ is the heat kernel of $\Delta ^{\alpha/2}$ on $(0,L)$ with Dirichlet boundary condition $u(t,x)=0,\ x\in (0,L)^{c}$.

For $p\geq 2$, Xie \cite{Xiebin} defined the $p$th energy of the solution at time $t$ by
\begin{align*}
\Phi_{p}(t,\lambda)=(\mathbb{E}\|u_{\lambda}(t)\|^{p}_{L^{p}})^{1/p}, \ t>0.
\end{align*}
\textbf{Definition 1.3.} The \emph{excitation index} of $u$ at time $t$ is given by
\begin{align*}
e_{p}(t):=\lim_{\lambda\rightarrow \infty}\frac{\log\log \Phi_{p}(t,\lambda)}{\log \lambda}.
\end{align*}

The parameter $\lambda>0$ in (\ref{sfh}) is called the level of noise (or noise level, for short).  Recently, effects of the noise level on nonlinear stochastic heat equations has attracted a great deal of research interest.
Khoshnevisan and Kim \cite{Kunwoo} studied a stochastic heat equation of the form
\begin{align*}
\frac{\partial}{\partial t}u=\mathcal{L}u+\lambda\sigma(u)\xi,
\end{align*}
where $\xi$ denotes space-time white noise on $\mathbb{R}_{+}\times G$, $\mathcal{L}$ is the generator of a L$\acute{e}$vy process on a locally compact Hausdorff Abelian group $G$, $\sigma: \mathbb{R}\rightarrow \mathbb{R}$ is Lipschitz
continuous, $\lambda$ is a large parameter. They showed that if $u$ is intermittent, the energy of the solution behaves generically as $\exp[\mbox{const}\cdot\lambda^{2}]$ when
$G$ is discrete and $\geq \exp[\mbox{const}\cdot\lambda^{4}]$ when $G$ is connected. Khoshnevisan and Kim \cite{KK}  studied the semilinear heat equation
\begin{align}\label{DK}
\partial_{t}u(t,x)=\partial_{xx}u(t,x)+\lambda \sigma(u(t,x))\xi
\end{align}
on the interval $[0,L]$, where $\xi$ is space-time white noise, $\sigma: \mathbb{R}\rightarrow \mathbb{R}$ is Lipschitz continuous.
The authors proved when the solution of \emph{Eq.} (\ref{DK}) is intermittent, the $L^{2}$-energy of the solution grows at lease $\exp(c\lambda^{2})$ and at most as
$\exp(c\lambda^{4})$ as $\lambda\rightarrow \infty$.  Foondun and Joseph \cite{FJ} considered \emph{Eq.} (\ref{DK}), they used Gaussian estimates for Dirichlet (Neumann) heat kernel
and renewal inequalities to show that the expected $L^{2}$-energy of the mild solution
is of order $\exp[\mbox{const}\cdot\lambda^{4}]$ as $\lambda\rightarrow \infty$. Foondun and Nualart \cite{FN} considered  a stochastic heat equation on an interval
\begin{align}\label{heat equation}
\partial _{t}u(t,x)=\frac{1}{2}\partial_{xx}u(t,x)+\lambda\sigma(u(t,x))\dot{w}(t,x),
\end{align}
where $\dot{w}$ denotes the space-time white noise.
They proved the second moment of the solution grows exponentially fast if the noise intensity $\lambda$ is large enough,
the second moment decays exponentially if $\lambda$ is small.
For $p\geq 2$, they proved the $p$th moments have a similar property. Xie \cite{Xiebin} studied \emph{Eq.} (\ref{heat equation}) on $[0,1]$, the author showed that for small noise level, the
$p$th moment of $\sup_{x\in [0,1]}|u(t,x)|$ is exponentially stable. For large noise level, the moment grows at least exponentially by an approach depending on the lower
bound of the global estimate for Dirichlet heat kernel. The $p$th energy of the solution at time $t$ is defined and the noise excitation index of the $p$th energy of
$u(t,x)$ is proved to be 4 as the noise level tends to infinity. Foondum, Tian and Liu \cite{FTL} studied nonlinear parabolic stochastic equations of the form
\begin{align*}
\partial_{t}u(t,x)=\mathfrak{L}u(t,x)+\lambda\sigma(u(t,x))\dot{w}(t,x)
\end{align*}
on the ball $B(0,\mathbb{R})$, where $\dot{w}$ denotes white noise on $(0,\infty)\times B(0,\mathbb{R})$, $\mathfrak{L}$ is the generator of an $\alpha$-stable process killed upon existing $B(0,\mathbb{R})$. The growth properties of the second moment of the solutions is obtained.
Their results extend those in \cite{FJ} and \cite{KK}.

In this paper, we prove that \emph{Eq.} (1.1) has a unique mild solution $u_{\lambda}(t,\cdot)$ in the Banach space $B_{p,\gamma}(\gamma<0)$, and then present an upper bound of the growth
rate of the mild solution for any $\lambda>0$. Our approach depends on the heat kernel estimates for the Dirichlet fractional Laplacian $\Delta^{\alpha/2}|_{(0,L)}:=-(-\Delta)^{\alpha/2}|_{(0,L)}$.
Assuming that $\lambda_{1}$ is the smallest eigenvalue of $(-\Delta)^{\alpha/2}$ on the interval $(0,L)$, for $p>2/(\alpha-1)$, $\beta\in (2/p,\alpha-1)$ and $\gamma\in (0,(2-\beta)\lambda_{1})$,
we prove that there exists $\lambda_{L}>0$ and we show that for all $\lambda\in (0,\lambda_{L})$ and $t\geq 0$,  the $p$th moment of $\|u_{\lambda}(t)\|_{L^{\infty}}$ is exponentially stable.
Then, we use the convolution-type inequalities and the asymptotic property of the Mittag-Leffler functions to prove
that the second moment of the mild solution grows faster than a Mittag-Leffler function. Then we give some properties of the $p$th moment and the $p$th energy of the mild solution and consider the non-linear noise excitability of \emph{Eq.} (1.1) for large noise level $\lambda$. At last, we show that
the noise excitation index of the solution $u_{\lambda}(t,x)$ with respect to $p$th energy $\Phi_{p}(t,\lambda)$ is $\frac{2\alpha}{\alpha-1}$.

The paper is organized as follows. In Section 2, we present some lemmas and preliminary facts, which will be used in the next sections.
In Section 3, we consider the existence and uniqueness of the mild solutions of \emph{Eq.} (1.1) and obtain some properties of the solutions. In Section 4, we prove that the $p$th
moment of the mild solution grows faster than a Mittag-Leffler function and obtain the excitation index of the mild solution of \emph{Eq.} (1.1).

Throughout this paper, we use $c_{0}$, $c_{1}$, $c_{2}$, $\ldots$ to denote generic constants, which may change from line to line.

\section{Preliminaries}
Let $p(t,x,y)$ be the heat kernel of $\Delta ^{\alpha/2}$ on $\mathbb{R}$.
We have the following inequality (see \cite{BGR})
\begin{align}\label{inequality}
0\leq p_{D}(t,x,y)\leq p(t,x,y) \ {\mbox{for all}} \ t>0, \ x,y \in \mathbb{R}.
\end{align}
For two nonnegative functions $f_{1}$ and $f_{2}$, the notion $f_{1}\asymp f_{2}$ means that $c_{1}f_{2}(x)\leq f_{1}(x)\leq c_{2}f_{2}(x)$, where $c_{1},c_{2}$ are positive constants.
It is well known that (see, e.g., \cite{TG,KT})
\begin{align*}
p(t,x,y)\asymp \big(t^{-1/\alpha}\wedge \frac{t}{|x-y|^{1+\alpha}}\big),
\end{align*}
that is, there exist constants $c_{1}$, $c_{2}$ such that for $t>0$, $x,y\in \mathbb{R}$,
\begin{align}\label{c1}
c_{1}\big(t^{-1/\alpha}\wedge \frac{t}{|x-y|^{1+\alpha}}\big)\leq p(t,x,y)\leq c_{2}\big(t^{-1/\alpha}\wedge \frac{t}{|x-y|^{1+\alpha}}\big),
\end{align}
where $c_{1}$ and $c_{2}$ are positive constants depending on $\alpha$. \\
Another equivalent form of (\ref{c1}) is
\begin{align}\label{alpha}
\frac{c_{1}t}{(t^{1/\alpha}+|y-x|)^{1+\alpha}}\leq p(t,x,y)\leq \frac{c_{2}t}{(t^{1/\alpha}+|y-x|)^{1+\alpha}},
\end{align}
where $t>0$, $x,y\in \mathbb{R}$, $c_{1}$ and $c_{2}$ are positive constants depending on $\alpha$. \\
\textbf{Lemma 2.1.} (Garsia-Rodemich-Rumsey Theorem, see \cite{FV}, Appendix A) Consider $f\in C([0,L],E)$, where $(E,d)$ is a complete metric space. Let $\Psi$ and $p$ be continuous strictly increasing functions on $[0,\infty)$ with $p(0)=\Psi(0)=0$ and $\Psi(x)\rightarrow \infty$ as $x\rightarrow \infty$. Then
\begin{align*}
\int_{0}^{L}\int_{0}^{L}\Psi\big(\frac{f(x)-f(y)}{p(|x-y|)}\big)dxdy\leq F,
\end{align*}
implies, for $0\leq x<y\leq L$,
\begin{align*}
d(f(x)-f(y))\leq 8\int_{0}^{y-x}\Psi^{-1}\big(\frac{4F}{u^{2}}\big)dp(u).
\end{align*}
In particular, if $\mbox{osc}(f,\delta)\equiv \sup\{d(f(x),f(y)): x,y\in [0,L], |x-y|\leq \delta\}$ denotes the modulus of continuity of $f$, we have
\begin{align*}
\mbox{osc}(f,\delta)\leq 8\int_{0}^{\delta}\Psi^{-1}\big(\frac{4F}{u^{2}}\big)dp(u).
\end{align*}
\textbf{Lemma 2.2.} Suppose that $\{u(x)\}_{x\in [0,L]}$ is a real valued stochastic process. If there exist $p\geq 1$ and positive constant $K,\delta$ such that
\begin{align}\label{Kolmogorov}
\mathbb{E}[|u(x)-u(y)|^{p}] \leq K|x-y|^{1+\delta},
\end{align}
then $\{u(x)\}_{x\in [0,L]}$ has a continuous modification, which will be still denoted by  $\{u(x)\}_{x\in [0,L]}$. For each $\varepsilon\in (0,\min\{\delta,1\})$, there exists a positive constant $\kappa$ depending only on $\delta$, $\varepsilon$ such that
\begin{align}\label{ukxy}
|u(x)-u(y)|\leq \kappa (4B)^{1/p}|x-y|^{(\delta-\varepsilon)/p},
\end{align}
where $B=B(p,\varepsilon,\delta)$ is a positive random variable defined by
\begin{align}\label{random}
B=\int_{0}^{L}\int_{0}^{L}\frac{|u(x)-u(y)|^{p}}{|x-y|^{2+\delta-\varepsilon}}dxdy.
\end{align}
In particular, the stochastic process $\{u(x)\}_{x\in [0,L]}$ has a $(\delta-\varepsilon)/p$-H$\ddot{\mbox{o}}$lder continuous modification. \\
\textbf{Proof.} By Lemma 2.1, we can finish the proof, which is similar to that of  Lemma 2.1 in \cite{Xiebin}, so it is omitted. \ \ $\Box$\\
\textbf{Remark 2.3.} There is an inaccuracy in the formula corresponding to (\ref{ukxy}) in \cite{Xiebin} (see \cite{Xiebin}, formula (2.5)).  $B$ should be changed to $4B$. The reason is that the  Theorem A.1  (Garsia, Rodemich and Rumsey inequality) in \cite{Xiebin} is not completely the same as Lemma 1.1 in \cite{GRR}. \\
\textbf{Lemma 2.4.} Suppose $b\geq 0, \beta> 0$ and $a(t)$ is a nonnegative function locally on $0\leq t<T$ (some $T\leq \infty$), and suppose $v(t)$ is nonnegative and locally integrable on $0\leq t<T$ with
\begin{align*}
v(t)\geq  a(t)+b\int_{0}^{t}(t-s)^{\beta-1}v(s)ds,
\end{align*}
on this interval; then
\begin{align*}
v(t)\geq a(t)+\theta\int_{0}^{t}F'_{\beta}(\theta(t-s))a(s)ds,\ 0\leq t< T,
\end{align*}
where
\begin{align*}
\theta=(b\Gamma(\beta))^{1/\beta}, \ F_{\beta}(z)=\sum_{n=0}^{\infty}\frac{z^{n\beta}}{\Gamma(n\beta+1)},\ F'_{\beta}(z)=\frac{d}{dz}F_{\beta}(z),
\end{align*}
If $a(t)\equiv a$, constant, then $v(t)\geq aF_{\beta}(\theta t)$. \\
\textbf{Proof.} The proof is similar to that of Lemma 7.1.1 in \cite{DH}, for the sake of completeness, we provide it in the following. Let
$B\phi(t)=b\int_{0}^{t}(t-s)^{\beta-1}\phi(s)ds$, $t\geq 0$, for locally integrable functions $\phi$. Then $v\geq a+Bv$ implies that
\begin{align*}
v\geq \sum_{k=0}^{n-1}B^{k}a+B^{n}v
\end{align*}
and
\begin{align*}
B^{n}v(t)=\frac{1}{\Gamma(n\beta)}\int_{0}^{t}(b\Gamma(\beta))^{n}(t-s)^{n\beta-1}u(s)ds\rightarrow 0
\end{align*}
as $n\rightarrow \infty$. Thus
\begin{align*}
v(t)&\geq a(t)+\int_{0}^{t}\{\sum_{n=1}^{\infty}\frac{1}{\Gamma(n\beta)}(b\Gamma(\beta))^{n}(t-s)^{n\beta-1}\}a(s)ds\\
&=a(t)+\theta\int_{0}^{t}F'_{\beta}(\theta(t-s))a(s)ds
\end{align*}
If $a(t)=a$ is a constant, then
\begin{align*}
v(t)&\geq a+a\theta\int_{0}^{t} F'_{\beta}(\theta s)ds\\
&\geq a+a\int_{0}^{\theta t}F'_{\beta}(s)ds\\
&=aF_{\beta}(\theta t).
\end{align*}
The proof is complete. \ $\Box$

\section{Properties of mild solutions to \emph{Eq.}(\ref{sfh}) in different Banach spaces}
In this section, we present the existence and uniqueness result of \emph{Eq.} (\ref{sfh}) in the Banach space $B_{p,\gamma}$. We prove that there exists $\gamma_{0}<0$,
for any $\gamma<\gamma_{0}$, \emph{Eq.} (\ref{sfh}) has a unique mild solution $u_{\lambda}(t,\cdot)\in B_{p,\gamma}$, and then we present an upper bound of the growth
rate of the solution for any $\lambda>0$. When $\lambda$ is small, under some assumptions, we show that there exists a constant $\lambda_{L}$ such that for $\lambda\in (0,\lambda_{L})$,
the $p$th moment of $\|u_{\lambda}(t)\|_{L^{\infty}}$ is exponential stable. \\
\textbf{Lemma 3.1.} Let $\beta\in (0,\alpha-1)$. There exists a positive constant $c(\alpha, \beta)$ that only depends on $\alpha$ and $\beta$ such that
\begin{align*}
\sup_{t\geq 0, \ x\in [0,L]}\int_{0}^{t}e^{\gamma s}s^{-2\beta/\alpha}\int_{0}^{L}|p_{D}(s,x,y)|^{2-\beta}dyds\leq c(\alpha, \beta)|\gamma|^{(\beta+1-\alpha)/\alpha} \ {\mbox{for any}} \ \gamma<0.
\end{align*}
\textbf{Proof.} By (\ref{inequality}), (\ref{alpha}), we have
\begin{align*}
&\int_{0}^{t}e^{\gamma s}s^{-2\beta/\alpha}\int_{0}^{L}|p_{D}(s,x,y)|^{2-\beta}dyds \\
&\leq \int_{0}^{t}e^{\gamma s}s^{-2\beta/\alpha}\int_{0}^{L}|p(s,x,y)|^{2-\beta}dyds\\
&\leq c_{2}^{2-\beta}(\alpha)\int_{0}^{t}e^{\gamma s}s^{-2\beta/\alpha}\int_{0}^{L}\big|\frac{s}{(s^{1/\alpha}+|y-x|)^{1+\alpha}}\big|^{2-\beta}dyds\\
&\leq c_{2}^{2-\beta}(\alpha)\int_{0}^{t}e^{\gamma s}s^{-2\beta/\alpha}s^{-(2-\beta)/\alpha}\int_{0}^{L}\frac{1}{(1+s^{-1/\alpha}|y-x|)^{(1+\alpha)(2-\beta)}}dyds\\
&\leq c_{2}^{2-\beta}(\alpha)\int_{0}^{t}e^{\gamma s}s^{-2\beta/\alpha}s^{-(2-\beta)/\alpha}\int_{\mathbb{R}}\frac{1}{(1+s^{-1/\alpha}|y-x|)^{(1+\alpha)(2-\beta)}}dyds\\
&\leq c_{2}^{2-\beta}(\alpha)\int_{0}^{t}e^{\gamma s}s^{-2\beta/\alpha}s^{-(2-\beta)/\alpha}s^{1/\alpha}\int_{\mathbb{R}}\frac{1}{(1+|\tau|)^{(1+\alpha)(2-\beta)}}d\tau ds\\
&\leq c_{2}^{2-\beta}(\alpha)\int_{0}^{t}e^{\gamma s}s^{-(\beta+1)/\alpha}\int_{\mathbb{R}}\frac{1}{(1+|\tau|)^{2}}d\tau ds\\
&\leq c_{2}^{2-\beta}(\alpha)\pi\int_{0}^{t}e^{\gamma s}s^{-(\beta+1)/\alpha}ds\\
&\leq c_{2}^{2-\beta}(\alpha)\pi |\gamma|^{(\beta+1-\alpha)/\alpha}\int_{0}^{\infty}e^{-s}s^{-(\beta+1)/\alpha}ds\\
&\leq c_{2}^{2-\beta}(\alpha)\pi\Gamma(\frac{\alpha-\beta-1}{\alpha})|\gamma|^{(\beta+1-\alpha)/\alpha},
\end{align*}
where $\Gamma(\cdot)$ is the Gamma function.  \ \ $\Box$ \

For any $u\in B_{p,\gamma}$, we define the stochastic convolution integral
\begin{align*}
Su(t,x)=\int_{0}^{t}\int_{0}^{L}p_{D}(t-s,x,y)u(s,y)w(dsdy).
\end{align*}
\textbf{Lemma 3.2.} Suppose $p>2/(\alpha-1)$ and $\gamma<0$. Then for each $\beta\in (2/p,\alpha-1)$, there exists a positive constant $c=c(\alpha,\beta,p)$ depending only on $\alpha$, $\beta$ and $p$, such that for all $u\in B_{p,\gamma}$,
\begin{align*}
\|Su\|^{p}_{p,\gamma}\leq c\|u\|^{p}_{p,\gamma}(|\gamma|^{(\beta+1-\alpha)p/2\alpha}+|\gamma|^{(1-\alpha)p/2\alpha}).
\end{align*}
\textbf{Proof.} By Burkh$\ddot{\mbox{o}}$lder's inequality, for any $p\in (0,\infty)$, there exists a positive constant $c(p)$ such that for any $x,y\in [0,L]$ and $t\geq 0$,
\begin{align*}
&\mathbb{E}|Su(t,x)-Su(t,y)|^{p}\\
&\leq c(p)\mathbb{E}\big(\int_{0}^{t}\int_{0}^{L}(p_{D}(t-s,x,z)-p_{D}(t-s,y,z))^{2}u^{2}(s,z)dzds\big)^{p/2},
\end{align*}
where $c(p)$ is a positive constant. By Minkowski's integral inequality (see \cite{Stein}, Appendice A.1), we have for $p\geq 2$, $\beta\in (0,2)$,
\begin{align}\label{Minkowski}
&\mathbb{E}|Su(t,x)-Su(t,y)|^{p}\nonumber\\
&\leq c(p)\big[\int_{0}^{t}\int_{0}^{L}(p_{D}(t-s,x,z)-p_{D}(t-s,y,z))^{2}(\mathbb{E}|u(s,z)|^{p})^{2/p}dzds\big]^{p/2}\nonumber\\
&=c(p)\big[\int_{0}^{t}\int_{0}^{L}|p_{D}(t-s,x,z)-p_{D}(t-s,y,z)|^{\beta}|p_{D}(t-s,x,z)\nonumber\\
&\quad-p_{D}(t-s,y,z)|^{2-\beta}\times(\mathbb{E}|u(s,z)|^{p})^{2/p}dzds\big]^{p/2}.
\end{align}
Note that for any $x,y\in \mathbb{R}$ and $f\in C^{1}(\mathbb{R})$, we have
\begin{align}\label{Taylor}
\int_{0}^{1}f'(x+r(y-x))(y-x)dr=f(y)-f(x).
\end{align}
By (\ref{Minkowski}), (\ref{Taylor}), we get
\begin{align}\label{kt}
&\mathbb{E}|Su(t,x)-Su(t,y)|^{p}\nonumber\\
&\leq c(p)|x-y|^{\beta p/2}\big[\int_{0}^{t}\int_{0}^{L}|\int_{0}^{1}\partial _{x}p_{D}(t-s,x+r(y-x),z)dr|^{\beta}|p_{D}(t-s,x,z)\nonumber\\
&\quad-p_{D}(t-s,y,z)|^{2-\beta}(\mathbb{E}|u(s,z)|^{p})^{2/p}dzds\big]^{p/2}\nonumber\\
&:=K_{1}(t,x,y)|x-y|^{\beta p/2},
\end{align}
where
\begin{align}\label{denot}
K_{1}(t,x,y)&=c(p)\big[\int_{0}^{t}\int_{0}^{L}|\int_{0}^{1}\partial _{x}p_{D}(t-s,x+r(y-x),z)dr|^{\beta}|p_{D}(t-s,x,z)\nonumber\\
&\quad-p_{D}(t-s,y,z)|^{2-\beta}\times(\mathbb{E}|u(s,z)|^{p})^{2/p}dzds\big]^{p/2}.
\end{align}
From the proof of Theorem 5.1 in \cite{CMN}, it follows that  for any $t>0$, $x,y\in(0,L)$,
\begin{align}\label{derivative}
|\partial_{x}p_{D}(t,x,y)|\leq c_{3}(\alpha)t^{-1/\alpha}p(t,x,y).
\end{align}
Then, (\ref{derivative}) together with (\ref{alpha}) yields that
\begin{align}\label{t}
|\partial_{x}p_{D}(t,x,y)|\leq c_{4}(\alpha)t^{-2/\alpha},
\end{align}
where $c_{4}(\alpha)$ is a positive constant depending only on $\alpha$. \\
By (\ref{denot}), (\ref{t}), we obtain
\begin{align}\label{convolution}
&K_{1}(t,x,y)\nonumber\\
&\leq c^{\beta}_{4}(\alpha)c(p)\big[\int_{0}^{t}\int_{0}^{L}(t-s)^{-2\beta/\alpha}|p_{D}(t-s,x,z)-p_{D}(t-s,y,z)|^{2-\beta}\nonumber\\
&\quad\cdot(\mathbb{E}|u(s,z)|^{p})^{2/p}dzds\big]^{p/2}\nonumber\\
&\leq c^{\beta}_{4}(\alpha)c(p)\big[\int_{0}^{t}\int_{0}^{L}(t-s)^{-2\beta/\alpha}|p_{D}(t-s,x,z)-p_{D}(t-s,y,z)|^{2-\beta}\nonumber\\
&\quad\cdot(\mathbb{E}||u(s)||^{p}_{L^{\infty}})^{2/p}dzds\big]^{p/2}\nonumber\\
&\leq c^{\beta}_{4}(\alpha)c(p)\big[\int_{0}^{t}\int_{0}^{L}(t-s)^{-2\beta/\alpha}|p_{D}(t-s,x,z)-p_{D}(t-s,y,z)|^{2-\beta}\nonumber\\
&\quad\cdot(e^{-\gamma s}e^{\gamma s}\mathbb{E}||u(s)||^{p}_{L^{\infty}})^{2/p}dzds\big]^{p/2}\nonumber\\
&\leq c^{\beta}_{4}(\alpha)c(p)\|u\|^{p}_{p,\gamma}\big[\int_{0}^{t}\int_{0}^{L}e^{-2\gamma s/p}(t-s)^{-2\beta/\alpha}|p_{D}(t-s,x,z)\nonumber\\
&\quad-p_{D}(t-s,y,z)|^{2-\beta}dzds\big]^{p/2}\nonumber\\
&\leq c^{\beta}_{4}(\alpha)c(p)e^{-\gamma t}\|u\|^{p}_{p,\gamma}\big[\int_{0}^{t}\int_{0}^{L}e^{2\gamma (t-s)/p}(t-s)^{-2\beta/\alpha}|p_{D}(t-s,x,z)\nonumber\\
&\quad-p_{D}(t-s,y,z)|^{2-\beta}dzds\big]^{p/2}\nonumber\\
&=c^{\beta}_{4}(\alpha)c(p)e^{-\gamma t}\|u\|^{p}_{p,\gamma}\big[\int_{0}^{t}\int_{0}^{L}e^{2\gamma s/p}s^{-2\beta/\alpha}|p_{D}(s,x,z)-p_{D}(s,y,z)|^{2-\beta}dzds\big]^{p/2}.
\end{align}
From Lemma 3.1 and (\ref{convolution}), it follows that
\begin{align}\label{constant}
K_{1}(t,x,y)&\leq c^{\beta}_{4}(\alpha)c(p)c^{p/2}(\alpha,\beta)e^{-\gamma t}\|u\|^{p}_{p,\gamma}|2\gamma/p|^{\frac{(\beta+1-\alpha)p}{2\alpha}}\nonumber\\
&=c_{5}(\alpha,\beta,p)e^{-\gamma t}\|u\|^{p}_{p,\gamma}|\gamma|^{(\beta+1-\alpha)p/2\alpha}\nonumber\\
&:=K_{2}(t).
\end{align}
where $c_{5}(\alpha,\beta,p)$ is a positive constant depending only on $\alpha,\beta, p$. \\
By (\ref{kt}) and (\ref{constant}), we have
\begin{align}\label{regular}
\mathbb{E}|Su(t,x)-Su(t,y)|^{p}&\leq c_{5}(\alpha,\beta,p)e^{-\gamma t}\|u\|^{p}_{p,\gamma}|\gamma|^{(\beta+1-\alpha)p/2\alpha}|x-y|^{\beta p/2}\nonumber\\
&=K_{2}(t)|x-y|^{\beta p/2}.
\end{align}
By Lemma 2.2, we have
\begin{align}\label{bound}
|Su(t,x)-Su(t,y)|^{p}\leq \kappa(4B(t))^{1/p}|x-y|^{\beta/2-(1+\varepsilon)/p}
\end{align}
where $\kappa$ is same as that in (\ref{estimate}),
\begin{align*}
B(t)=\int_{0}^{L}\int_{0}^{L}\frac{|Su(t,x)-Su(t,y)|^{p}}{|x-y|^{1+\beta p/2-\varepsilon}}dxdy,
\end{align*}
and $\varepsilon\in (0, \min\{\beta p/2-1,1\})$. \\
Fix $y=y_{0}\in [0,L]$. By (\ref{bound}), we get
\begin{align*}
|Su(t,x)|\leq \kappa((4B(t))^{1/p}|x-y_{0}|^{\beta/2-(1+\varepsilon)/p}+|Su(t,y_{0})|,
\end{align*}
then for $x,y_{0}\in [0,L]$,
\begin{align}\label{sup}
|Su(t,x)|\leq \kappa(4B(t))^{1/p}L^{\beta/2-(1+\varepsilon)/p}+|Su(t,y_{0})|,
\end{align}
where $\beta/2-(1+\varepsilon)/p>0$.
Taking the $p$th moments to both sides of (\ref{sup}) and consider the inequality $|a+b|^{p}\leq 2^{p-1}(|a|^{p}+|b|^{p})$ ($a,b\in \mathbb{R}$), we get
\begin{align}\label{analogous}
|Su(t,x)|^{p}\leq 2^{p-1}(4\kappa^{p}(B(t))L^{\beta p/2-(1+\varepsilon)}+|Su(t,y_{0})|^{p}).
\end{align}
Taking expectations of both sides of the above inequality, we have
\begin{align}\label{expectation}
\mathbb{E}\sup_{x\in [0,L]}|Su(t,x)|^{p}\leq 2^{p-1}\kappa^{p}L^{\beta p/2-(1+\varepsilon)}\mathbb{E}(B(t))+2^{p-1}\mathbb{E}|Su(t,y_{0})|^{p}.
\end{align}
By (\ref{regular}), we have
\begin{align}\label{norm}
\mathbb{E}(B(t))&\leq K_{2}(t)\int_{0}^{L}\int_{0}^{L}\frac{1}{|x-y|^{1-\varepsilon}}dxdy\nonumber\\
&\leq \frac{2L^{\varepsilon+1}}{\varepsilon(\varepsilon+1)}K_{2}(t)\nonumber\\
&=\frac{2L^{\varepsilon+1}}{\varepsilon(\varepsilon+1)}c_{5}(\alpha,\beta,p)e^{-\gamma t}\|u\|^{\beta}_{p,\gamma}|\gamma|^{(\beta+1-\alpha)p/2\alpha}.
\end{align}
By Burkh$\ddot{\mbox{o}}$lder's inequality, Minkowski's integral inequality and the semigroup property, we obtain
\begin{align}\label{ty}
\mathbb{E}|Su(t,y_{0})|^{p}&=\mathbb{E}\big|\int_{0}^{t}\int_{0}^{L}p_{D}(t-s,y_{0},y)u(s,y)w(dsdy)\big|^{p}\nonumber\\
&\leq c_{p}\mathbb{E}\big[\int_{0}^{t}\int_{0}^{L}p^{2}_{D}(t-s,y_{0},y)u^{2}(s,y)dyds\big]^{p/2}\nonumber\\
&\leq c_{p}\big[\int_{0}^{t}\int_{0}^{L}p^{2}_{D}(t-s,y_{0},y)(\mathbb{E}|u(s,y)|^{p})^{2/p}dyds\big]^{p/2}\nonumber\\
&\leq c_{p}\big[\int_{0}^{t}\int_{0}^{L}p^{2}_{D}(t-s,y_{0},y)e^{-2\gamma s/p}(e^{\gamma s}\mathbb{E}\|u(s)\|^{p}_{L^{\infty}})^{2/p}dyds\big]^{p/2}\nonumber\\
&=c_{p}\|u\|^{p}_{p,\gamma}\big[\int_{0}^{t}e^{-2\gamma s/p}\big(\int_{0}^{L}p^{2}_{D}(t-s,y_{0},y)dy\big)ds\big]^{p/2}\nonumber\\
&=c_{p}e^{-\gamma t}\|u\|^{p}_{p,\gamma}\big[\int_{0}^{t}e^{2\gamma (t-s)/p}p_{D}(2(t-s),y_{0},y_{0})ds\big]^{p/2}\nonumber\\
&=c_{p}e^{-\gamma t}\|u\|^{p}_{p,\gamma}\big[\int_{0}^{t}e^{2\gamma s/p}p_{D}(2s,y_{0},y_{0})ds\big]^{p/2}\nonumber\\
&\leq c_{p}c^{p/2}_{7}(\alpha)2^{-p/2\alpha}e^{-\gamma t}\|u\|^{p}_{p,\gamma}\big(\int_{0}^{t}e^{2\gamma s/p}s^{-1/\alpha}ds\big)^{p/2}\nonumber\\
&\leq c_{p}c^{p/2}_{7}(\alpha)2^{-p/2\alpha}e^{-\gamma t}\|u\|^{p}_{p,\gamma}\big((\frac{p}{2|\gamma|})^{1-1/\alpha}\int_{0}^{\infty}e^{-s}s^{-1/\alpha}ds \big)^{p/2}\nonumber\\
&\leq c_{p}c^{p/2}_{7}(\alpha)2^{-p/2\alpha}\big((\frac{p}{2})^{1-1/\alpha}\Gamma(1-\frac{1}{\alpha}) \big)^{p/2}e^{-\gamma t}\|u\|^{p}_{p,\gamma}|\gamma|^{(1-\alpha)p/2\alpha}.
\end{align}
Putting (\ref{norm}) and (\ref{ty}) into (\ref{expectation}), we have
\begin{align}\label{complete}
\mathbb{E}\sup_{x\in [0,L]}|Su(t,x)|^{p}\leq c_{8}(\alpha,\beta,p)e^{-\gamma t}\|u\|^{p}_{p,\gamma}(|\gamma|^{(\beta+1-\alpha)p/2\alpha}+|\gamma|^{(1-\alpha)p/2\alpha}),
\end{align}
where $c_{8}(\alpha,\beta,p)$ is a positive constant that only depends on $\alpha$, $\beta$ and $p$. From (\ref{complete}) and (\ref{gamma}), it follows that
\begin{align}\label{proof}
\|Su\|^{p}_{p,\gamma}\leq c_{8}(\alpha,\beta,p)\|u\|^{p}_{p,\gamma}(|\gamma|^{(\beta+1-\alpha)p/2\alpha}+|\gamma|^{(1-\alpha)p/2\alpha}).
\end{align}
The proof is complete. \  $\Box$ \\
\textbf{Lemma 3.3.} Suppose $p>2/(\alpha-1)$ and $\gamma<0$. Then for each $\beta\in (2/p,\alpha-1)$, there exists a positive constant $c=c(\alpha,\beta,p)$ depending only on $\alpha$, $\beta$ and $p$, such that for any $u,v\in B_{p,\gamma}$,
\begin{align*}
\|Su-Sv\|^{p}_{p,\gamma}\leq c\|u-v\|^{p}_{p,\gamma}(|\gamma|^{(\beta+1-\alpha)p/2\alpha}+|\gamma|^{(1-\alpha)p/2\alpha}).
\end{align*}
\textbf{Proof.} The proof is similar to that of Lemma 3.2, we omit it. \ \ $\Box$ \ \\
\textbf{Theorem 3.4.}  Suppose $p>2/(\alpha-1)$. Then exists $\gamma_{0}<0$ such that for all $\gamma<\gamma_{0}$, $Eq$. (\ref{sfh}) has a unique mild solution $u_{\lambda}(t,\cdot)\in B_{p,\gamma}$. For any $\lambda>0$ and all $\gamma<\gamma_{0}$, the following inequality holds
\begin{align*}
\lim_{t\rightarrow \infty}\sup\frac{1}{t}\log\mathbb{E}\|u_{\lambda}(t)\|^{p}_{L^{\infty}}\leq -\gamma.
\end{align*}
That is, the growth of $u_{\lambda}(t,x)$ in time $t$ is at most in an exponential rate in the $p$-moment sense. \\
\textbf{Proof.} Define the operator $T$ as follows
\begin{align*}
Tu(t,x)=\int_{0}^{L}p_{D}(t,x,y)u_{0}(y)dy+\int_{0}^{t}\int_{0}^{L}p_{D}(t-s,x,y)\lambda \sigma(u(s,y))w(dsdy),
\end{align*}
where $u\in B_{p,\gamma}$. \\
Since $u_{0}$ is non-random, for any $\gamma<0$, we have
\begin{align*}
&\big\|\int_{0}^{L}p_{D}(t,x,y)u_{0}(y)dy\big\|^{p}_{p,\gamma}\nonumber\\
&\quad\leq \|u_{0}\|^{p}_{L^{\infty}}\sup_{t\geq 0}e^{\gamma t}\big|\int_{0}^{L}p_{D}(t,x,y)dy\big|^{p}\nonumber\\
&\quad\leq \|u_{0}\|^{p}_{L^{\infty}}\sup_{t\geq 0}e^{\gamma t}\big|\int_{0}^{L}p(t,x,y)dy\big|^{p}\nonumber\\
&\quad\leq \|u_{0}\|^{p}_{L^{\infty}}.
\end{align*}
Define the stochastic convolution
\begin{align*}
S\sigma(u(t,x))=\int_{0}^{t}\int_{0}^{L}p_{D}(t-s,x,y)\lambda \sigma(u(s,y))w(dsdy).
\end{align*}
By Assumption 1.2, $\sigma(u)\leq L_{\sigma}|u|$. From Lemma 3.2, it follows that for any $u\in B_{p,\gamma}$,
\begin{align*}
\|Tu\|^{p}_{p,\gamma}\leq 2^{p-1}\|u_{0}\|^{p}_{L^{\infty}}+2^{p-1}\lambda^{p}c(\alpha,\beta,p,\kappa,L)L^{p}_{\sigma}\|u\|^{p}_{p,\gamma}(|\gamma|^{(\beta+1-\alpha)p/2\alpha}+|\gamma|^{(1-\alpha)p/2\alpha}),
\end{align*}
where  $p>2/(\alpha-1)$, $\beta\in (2/p,\alpha-1)$ and $\gamma<0$. Then we see that $T$ maps $B_{p,\gamma}$ into $B_{p,\gamma}$ for any $\gamma<0$.
For any $u,v\in B_{p,\gamma}$, by Lemma 3.3, we have
\begin{align}\label{contract}
\|Tu-Tv\|^{p}_{p,\gamma}&=L^{p}_{\sigma}\|S\sigma(u)-S\sigma(v)\|^{p}_{p,\gamma}\nonumber\\
&\leq c(\alpha,\beta,p)\lambda^{p}L^{p}_{\sigma}(|\gamma|^{(\beta+1-\alpha)p/2\alpha}+|\gamma|^{(1-\alpha)p/2\alpha})\|u-v\|^{p}_{p,\gamma},
\end{align}
Since $2/p+1<\alpha<2$, $0<\beta<\alpha-1$, we can choose
a $\gamma_{0}<0$ such that for any $\gamma<\gamma_{0}$,
\begin{align*}
0<c(\alpha,\beta,p)\lambda^{p}L^{p}_{\sigma}(|\gamma|^{(\beta+1-\alpha)p/2\alpha}+|\gamma|^{(1-\alpha)p/2\alpha})<1.
\end{align*}
From the contraction mapping theorem, it follows  that $T$ has a unique fixed point on $B_{p,\gamma}$, therefore \emph{Eq.}(\ref{sfh}) has a unique mild solution $u_{\lambda}(t,\cdot)\in B_{p,\gamma}$.
Since $\|u_{\lambda}\|_{p,\gamma}=\|Tu_{\lambda}\|_{p,\gamma}<\infty$, it is easy to see that $\lim_{t\rightarrow \infty}\sup\frac{1}{t}\log\mathbb{E}\|u_{\lambda}(t)\|^{p}_{L^{\infty}}\leq -\gamma$.  The proof is complete. \ \ $\Box$ \ \\
\textbf{Lemma 3.5} (see \cite{Chen}), Theorem 1.1.) Let $D$ be a $C^{1,1}$ open subset of $\mathbb{R}^{d}$ with $d\geq 1$ and $\delta_{D}(x)$ the Euclidean distance between $x$ and $D^{c}$. \\
(i) For every $T>0$, on $(0,T]\times D\times D$,
\begin{align*}
p_{D}(t,x,y)\asymp \big(1\wedge \frac{\delta_{D}(x)^{\alpha/2}}{\sqrt{t}}\big)\big(1\wedge \frac{\delta_{D}(y)^{\alpha/2}}{\sqrt{t}}\big)\big(t^{-d/\alpha}\wedge \frac{t}{|x-y|^{d+\alpha}}\big).
\end{align*}
(ii) Suppose in addition that $D$ is bounded. For every $T>0$, there are positive constants $c_{1}< c_{2}$ such that on $[T,\infty)\times D\times D$,
\begin{align*}
c_{1}e^{-\lambda_{0}t}\delta_{D}(x)^{\alpha/2}\delta_{D}(y)^{\alpha/2}\leq p_{D}(t,x,y)\leq c_{2}e^{-\lambda_{0}t}\delta_{D}(x)^{\alpha/2}\delta_{D}(y)^{\alpha/2},
\end{align*}
where $\lambda_{0}>0$ is the smallest eigenvalue of the Dirichlet fractional Laplacian $(-\Delta)^{\alpha/2}|_{D}$, $p_{D}(t,x,y)$ is the Dirichlet fractional heat kernel.
For the definition of $C^{1,1}$ open set, we refer to  \cite{Chen}.

Throughout this paper,  we assume that  $\lambda_{1}$ is the smallest eigenvalue of the Dirichlet fractional Laplacian $(-\Delta)^{\alpha/2}|_{(0,L)}$. \\
\textbf{Lemma 3.6} Suppose $\beta\in (0,\alpha-1)$ and $\gamma\in (0,(2-\beta)\lambda_{1})$. Then there exists a positive constant $c(\alpha,\beta)$, such that for any $t\geq 0$ and $x\in (0,L)$,
\begin{align*}
\int_{0}^{t}e^{\gamma s}s^{-2\beta/\alpha}\int_{0}^{L}|p_{D}(s,x,y)|^{2-\beta}dyds \leq c(\alpha,\beta) \big(\gamma^{(\beta+1-\alpha)/\alpha}+\frac{1}{(2-\beta)\lambda_{1}-\gamma}\big),
\end{align*}
\textbf{Proof.} By (\ref{inequality}), (\ref{alpha}), similar to the proof of Lemma 2.3,  we have
\begin{align}\label{T}
&\int_{0}^{1}e^{\gamma s}s^{-2\beta/\alpha}\int_{0}^{L}|p_{D}(s,x,y)|^{2-\beta}dyds\nonumber\\
&\leq \int_{0}^{1}e^{\gamma s}s^{-2\beta/\alpha}\int_{0}^{L}|p(s,x,y)|^{2-\beta}dyds\nonumber\\
&\leq c_{2}^{2-\beta}(\alpha)\pi\int_{0}^{1}e^{\gamma s}s^{-(\beta+1)/\alpha}ds\nonumber\\
&=c_{2}^{2-\beta}(\alpha)\pi\gamma^{(\beta+1-\alpha)/\alpha}\int_{0}^{2\lambda_{1}}e^{s}s^{-(\beta+1)/\alpha}ds.
\end{align}
Note that $-(\beta+1)/\alpha\in (-1,0)$, then
\begin{align}\label{first bound}
\int_{0}^{T}e^{\gamma s}s^{-2\beta/\alpha}\int_{0}^{L}|p_{D}(s,x,y)|^{2-\beta}dyds \leq c_{1}(\alpha,\beta)\gamma^{(\beta+1-\alpha)/\alpha}.
\end{align}
By Lemma 3.5, for $\gamma\in (2-\beta)\lambda_{1}$, we have
\begin{align}\label{second bound}
&\int_{1}^{\infty}e^{\gamma s}s^{-2\beta/\alpha}\int_{0}^{L}|p_{D}(s,x,y)|^{2-\beta}dyds\nonumber\\
&\leq c^{2-\beta}_{2}\int_{1}^{\infty}e^{\gamma s}s^{-2\beta/\alpha}\int_{0}^{L}e^{-(2-\beta)\lambda_{1}s}\delta_{D}(x)^{\alpha(2-\beta)/2}\delta_{D}(y)^{\alpha(2-\beta)/2}dyds\nonumber\\
&\leq  c^{2-\beta}_{2}L^{\alpha(2-\beta)+1}\int_{1}^{\infty}e^{(\gamma-(2-\beta)\lambda_{1})s}s^{-2\beta/\alpha}ds\nonumber\\
&\leq c^{2-\beta}_{2}L^{\alpha(2-\beta)+1}\int_{1}^{\infty}e^{(\gamma-(2-\beta)\lambda_{1})s}ds\nonumber\\
&\leq c_{2}(\alpha,\beta)/((2-\beta)\lambda_{1}-\gamma).
\end{align}
By (\ref{first bound}) and (\ref{second bound}), the proof is complete. \ \  $\Box$ \ \\
\textbf{Theorem 3.7.} Suppose $p>2/(\alpha-1)$. There exists $\lambda_{L}>0$ such that for all $\lambda\in (0,\lambda_{L})$,
\begin{align}\label{lu}
-\infty<\lim_{t\rightarrow \infty}\sup\frac{1}{t}\log\mathbb{E}|u_{\lambda}(t,x)|^{p}\leq \lim_{t\rightarrow \infty}\sup\frac{1}{t}\log\mathbb{E}\|u_{\lambda}(t)\|^{p}_{L^{\infty}}<0.
\end{align}
\textbf{Proof.} Define
\begin{align}\label{f}
u_{\lambda}(t,x)=\int_{0}^{L}p_{D}(t,x,y)u_{0}(y)dy+\int_{0}^{t}\int_{0}^{L}p_{D}(t-s,x,y)\lambda \sigma(u(s,y))w(dsdy).
\end{align}
Taking the second moments to both sides of (\ref{f}), we get
\begin{align}\label{than}
\mathbb{E}|u_{\lambda}(t,x)|^{2}=&\big|\int_{0}^{L}p_{D}(t,x,y)u_{0}(y)dy\big|^{2}\nonumber\\
&\quad+\lambda^{2}\int_{0}^{t}\int_{0}^{L}p_{D}(t-s,x,y)\mathbb{E}|\sigma(u(s,y))|^{2}dyds.
\end{align}
By (ii) of Lemma 2.6, we see that $\big|\int_{0}^{L}p_{D}(t,x,y)u_{0}(y)dy\big|^{2}$ decays exponentially fast with time.
By (\ref{than}), $\mathbb{E}|u_{\lambda}(t,x)|^{2}$ can not decay faster than exponential.
Since $p>2/(\alpha-1)$, by Jensen's inequality, we have
\begin{align}\label{compare}
\mathbb{E}|u_{\lambda}(t,x)|^{p}=\mathbb{E}(|u_{\lambda}(t,x)|^{2})^{p/2}\geq (\mathbb{E}|u_{\lambda}(t,x)|^{2})^{p/2},
\end{align}
then
\begin{align*}
\lim_{t\rightarrow \infty}\sup\frac{1}{t}\log\mathbb{E}|u_{\lambda}(t,x)|^{p}>-\infty.
\end{align*}
By (\ref{gamma}), to prove (\ref{lu}), we only to show that for some $\gamma\in (0,(2-\beta)\lambda_{1})$, $\beta\in(2/p,\alpha-1)$, there exists $\lambda_{L}>0$ such that for any $\lambda\in (0,\lambda_{L})$, $u_{\lambda}\in B_{p,\gamma}$. \\
For any $\gamma\in (0,\lambda_{1}p)$, by (\ref{inequality}) and (ii) of Lemma 3.5, we have
\begin{align}\label{pdtxy}
&\big\|\int_{0}^{L}p_{D}(t,x,y)u_{0}(y)dy\big\|^{p}_{p,\gamma}\nonumber\\
&\quad\leq \|u_{0}\|^{p}_{L^{\infty}}\sup_{t\geq 0}e^{\gamma t}\big|\int_{0}^{L}p_{D}(t,x,y)u_{0}(y)dy\big|^{p}\nonumber\\
&\quad\leq \|u_{0}\|^{p}_{L^{\infty}}\big(\sup_{t\in[0,1]}e^{\gamma t}\big|\int_{0}^{L}p_{D}(t,x,y)dy\big|^{p}+\sup_{t\in(1,\infty)}e^{\gamma t}\big|\int_{0}^{L}p_{D}(t,x,y)dy\big|^{p}\big)\nonumber\\
&\quad\leq \|u_{0}\|^{p}_{L^{\infty}}\big(\sup_{t\in[0,1]}e^{\gamma t}\big|\int_{0}^{L}p(t,x,y)dy\big|^{p}+\sup_{t\in(1,\infty)}e^{\gamma t}\big|\int_{0}^{L}p_{D}(t,x,y)dy\big|^{p}\big)\nonumber\\
&\quad\leq \|u_{0}\|^{p}_{L^{\infty}}\big(e^{\gamma}+c_{2}\delta_{D}(x)^{\alpha p/2}\delta_{D}(y)^{\alpha p/2}\sup_{t\in(1,\infty)}e^{(\gamma-\lambda_{1}p)t}\big)\nonumber\\
&\quad\leq \|u_{0}\|^{p}_{L^{\infty}}\big(e^{\gamma}+c_{2}L^{\alpha p}\big),
\end{align}
which implies that for $\gamma\in (0,\lambda_{1}p)$, $\int_{0}^{L}p_{D}(t,x,y)u_{0}(y)dy\in B_{p,\gamma}$.\\
Define the mapping $S_{\lambda}$
\begin{align*}
S_{\lambda}(\sigma(u_{\lambda}(t,x)))=\int_{0}^{t}\int_{0}^{L}p_{D}(t-s,x,y)\lambda \sigma(u(s,y))w(dsdy).
\end{align*}
By Assumption 1.2, $\sigma(u)\leq L_{\sigma}|u|$. Similar to (\ref{kt}), we have
\begin{align}\label{su}
&\mathbb{E}|S_{\lambda}(\sigma(u_{\lambda}(t,x)))-S_{\lambda}(\sigma(u_{\lambda}(t,y)))|^{p}\nonumber\\
&\leq c(p)\lambda^{p}L^{p}_{\sigma}|x-y|^{\beta p/2}\big[\int_{0}^{t}\int_{0}^{L}|\int_{0}^{1}\partial _{x}p_{D}(t-s,x+r(y-x),z)dr|^{\beta}|p_{D}(t-s,x,z)\nonumber\\
&\quad-p_{D}(t-s,y,z)|^{2-\beta}\times(\mathbb{E}|u(s,z)|^{p})^{2/p}dzds\big]^{p/2}\nonumber\\
&:=K_{3}(t,x,y)|x-y|^{\beta p/2},
\end{align}
where
\begin{align}\label{denote}
K_{3}(t,x,y)&=c(p)\lambda^{p}L^{p}_{\sigma}\big[\int_{0}^{t}\int_{0}^{L}|\int_{0}^{1}\partial _{x}p_{D}(t-s,x+r(y-x),z)dr|^{\beta}\nonumber\\
&\quad\cdot|p_{D}(t-s,x,z)-p_{D}(t-s,y,z)|^{2-\beta}\times(\mathbb{E}|u(s,z)|^{p})^{2/p}dzds\big]^{p/2}.
\end{align}
Similar to (\ref{convolution}) and by Lemma 3.6, we have
\begin{align}\label{ktxy}
&K_{3}(t,x,y)\nonumber\\
&\leq c^{\beta}_{4}(\alpha)c(p)\lambda^{p}L^{p}_{\sigma}e^{-\gamma t}\|u\|^{p}_{p,\gamma}\big[\int_{0}^{t}\int_{0}^{L}e^{2\gamma s/p}s^{-2\beta/\alpha}|p_{D}(s,x,z)\nonumber\\
&\quad-p_{D}(s,y,z)|^{2-\beta}dzds\big]^{p/2}\nonumber\\
&\leq  c^{\beta}_{4}(\alpha)c(p)\lambda^{p}L^{p}_{\sigma}c^{p/2}(\alpha,\beta)e^{-\gamma t}\|u\|^{p}_{p,\gamma}\big((2\gamma/p)^{(\beta+1-\alpha)/\alpha}+\frac{1}{(2-\beta)\lambda_{1}-\gamma}\big)\nonumber\\
&:=K_{4}(t).
\end{align}
Fix $y=y_{0}\in [0,L]$, similar to (\ref{expectation}), we have
\begin{align}\label{expect}
\mathbb{E}\sup_{x\in [0,L]}|S_{\lambda}\sigma(u_{\lambda}(t,x))|^{p}\leq 2^{p-1}\kappa^{p}L^{\beta p/2-(1+\varepsilon)}\mathbb{E}(B(t))+2^{p-1}\mathbb{E}|S_{\lambda}\sigma(u_{\lambda}(t,y_{0}))|^{p}.
\end{align}
Similar to (\ref{norm}), by (\ref{ktxy}), we have
\begin{align}\label{bk}
\mathbb{E}(B(t))&\leq K_{4}(t)\int_{0}^{L}\int_{0}^{L}\frac{1}{|x-y|^{1-\varepsilon}}dxdy\nonumber\\
&\leq \frac{2L^{\varepsilon+1}}{\varepsilon(\varepsilon+1)}K_{4}(t)\nonumber\\
&=\frac{2L^{\varepsilon+1}}{\varepsilon(\varepsilon+1)}c^{\beta}_{4}(\alpha)c(p)\lambda^{p}L^{p}_{\sigma}c^{p/2}(\alpha,\beta)e^{-\gamma t}\|u\|^{p}_{p,\gamma}\big((2\gamma/p)^{(\beta+1-\alpha)/\alpha}\nonumber\\
&\quad+\frac{1}{(2-\beta)\lambda_{1}-\gamma}\big).
\end{align}
Similar to the proof of (\ref{ty}), we get
\begin{align}\label{sty}
&\mathbb{E}|S_{\lambda}\sigma(u_{\lambda}(t,y_{0}))|^{p}\nonumber\\
&=\mathbb{E}\big|\int_{0}^{t}\int_{0}^{L}p_{D}(t-s,y_{0},y)\lambda\sigma(u(s,y))w(dsdy)\big|^{p}\nonumber\\
&\leq c_{p}\mathbb{E}\big[\int_{0}^{t}\int_{0}^{L}p^{2}_{D}(t-s,y_{0},y)\lambda^{2}\sigma^{2}(u(s,y))dyds\big]^{p/2}\nonumber\\
&\leq c_{p}\lambda^{p}L^{p}_{\sigma}\mathbb{E}\big[\int_{0}^{t}\int_{0}^{L}p^{2}_{D}(t-s,y_{0},y)\lambda^{2}u^{2}(s,y)dyds\big]^{p/2}\nonumber\\
&\leq c_{p}\lambda^{p}L^{p}_{\sigma}\big[\int_{0}^{t}\int_{0}^{L}p^{2}_{D}(t-s,y_{0},y)(\mathbb{E}|u(s,y)|^{p})^{2/p}dyds\big]^{p/2}\nonumber\\
&\leq c_{p}\lambda^{p}L^{p}_{\sigma}\big[\int_{0}^{t}\int_{0}^{L}p^{2}_{D}(t-s,y_{0},y)e^{-2\gamma s/p}(e^{\gamma s}\mathbb{E}\|u(s)\|^{p}_{L^{\infty}})^{2/p}dyds\big]^{p/2}\nonumber\\
&=c_{p}\lambda^{p}L^{p}_{\sigma}\|u\|^{p}_{p,\gamma}\big[\int_{0}^{t}e^{-2\gamma s/p}\big(\int_{0}^{L}p^{2}_{D}(t-s,y_{0},y)dy\big)ds\big]^{p/2}\nonumber\\
&=c_{p}\lambda^{p}L^{p}_{\sigma}e^{-\gamma t}\|u\|^{p}_{p,\gamma}\big[\int_{0}^{t}e^{2\gamma (t-s)/p}p_{D}(2(t-s),y_{0},y_{0})ds\big]^{p/2}\nonumber\\
&=c_{p}\lambda^{p}L^{p}_{\sigma}e^{-\gamma t}\|u\|^{p}_{p,\gamma}\big[\int_{0}^{t}e^{2\gamma s/p}p_{D}(2s,y_{0},y_{0})ds\big]^{p/2}.
\end{align}
By (\ref{inequality}), (\ref{c1}), we have
\begin{align}\label{pd}
&\int_{0}^{1}e^{2\gamma s/p}p_{D}(2s,y_{0},y_{0})ds\nonumber\\
&\leq c_{2}(\alpha)2^{-1/\alpha}\int_{0}^{1}e^{2\gamma s/p}s^{-1/\alpha}ds\nonumber\\
&\leq c_{2}(\alpha)(p/2)^{(\alpha-1)/\alpha}\gamma^{(1-\alpha)/\alpha}\int_{0}^{2\gamma/p}e^{s}s^{-1/\alpha}ds\nonumber\\
&\leq c_{2}(\alpha)(p/2)^{(\alpha-1)/\alpha}\gamma^{(1-\alpha)/\alpha}\int_{0}^{2\lambda_{1}}e^{s}s^{-1/\alpha}ds.
\end{align}
By (ii) of Lemma 2.6, we get
\begin{align}\label{other}
&\int_{1}^{\infty}e^{2\gamma s/p}p_{D}(2s,y_{0},y_{0})ds\nonumber\\
&\leq c_{2}L^{\alpha}\int_{1}^{\infty}e^{(2\gamma/p-2\lambda_{1})s}ds\nonumber\\
&=\frac{c_{2}p}{2\lambda_{1}p-2\gamma}.
\end{align}
By (\ref{pd}) and (\ref{other}), we have
\begin{align}\label{pdy}
\int_{0}^{t}e^{2\gamma s/p}p_{D}(2s,y_{0},y_{0})ds\leq c(\alpha,p)\big(\gamma^{(1-\alpha)/\alpha}+\frac{1}{\lambda_{1}p-\gamma}\big),
\end{align}
where $c(\alpha,p,\lambda_{1})$ depends on $\alpha$, $p$ and $\lambda_{1}$. \\
Putting (\ref{pdy}) into (\ref{sty}), we obtain
\begin{align}\label{lambda}
\mathbb{E}|S_{\lambda}\sigma(u_{\lambda}(t,y_{0}))|^{p}\leq c_{p}\lambda^{p}L^{p}_{\sigma}c^{p/2}(\alpha,p)\big(\gamma^{(1-\alpha)/\alpha}+\frac{1}{\lambda_{1}p-\gamma}\big)^{p/2}e^{-\gamma t}\|u\|^{p}_{p,\gamma}.
\end{align}
Note that $2-\beta<2<p$, by (\ref{expect}), (\ref{bk}) and (\ref{lambda}), for $\gamma\in (0,(2-\beta)\lambda_{1})$, we have
\begin{align}\label{es}
&\sup_{x\in [0,L]}\mathbb{E}|S_{\lambda}\sigma(u_{\lambda}(t,x))|^{p}\nonumber\\
&\leq c_{1}e^{-\gamma t}\|u_{\lambda}\|_{p,\gamma}^{p}\big((2\gamma/p)^{(\beta+1-\alpha)/\alpha}+\frac{1}{(2-\beta)\lambda_{1}-\gamma}\big)\nonumber\\
&\quad+c_{3}\lambda^{p}e^{-\gamma t}\|u_{\lambda}\|_{p,\gamma}^{p}\big(\gamma^{(1-\alpha)/\alpha}+\frac{1}{\lambda_{1}p-\gamma}\big)^{p/2},
\end{align}
where $c_{1}=c(\alpha,\beta,p)$, $c_{3}=c(\alpha,p)$.
This implies that for $\gamma\in (0,(2-\beta)\lambda_{1})$, $S_{\lambda}\sigma(u_{\lambda}(t,y_{0}))\in B_{p,\gamma}$. \\
By (\ref{f}), (\ref{pdtxy}) and (\ref{es}), we have
\begin{align*}
\|u_{\lambda}\|_{p,\gamma}^{p}&\leq 2^{p-1}(e^{\gamma}+c_{2}L^{\alpha}p)\|u_{0}\|^{p}_{L^{\infty}}\nonumber\\
&\quad+2^{p-1}(c_{1}+c_{3}\lambda^{p}L_{\sigma}^{p})\|u_{\lambda}\|_{p,\gamma}^{p}\big[(2\gamma/p)^{(\beta+1-\alpha)/\alpha}+\frac{1}{(2-\beta)\lambda_{1}-\gamma}\big)\nonumber\\
&\quad+\big(\gamma^{(1-\alpha)/\alpha}+\frac{1}{\lambda_{1}p-\gamma}\big)^{p/2}\big]<\infty.
\end{align*}
We can choose $\lambda$  sufficiently small such that
\begin{align*}
&2^{p-1}(c_{1}+c_{3}\lambda^{p}L_{\sigma}^{p})\big[(2\gamma/p)^{(\beta+1-\alpha)/\alpha}+\frac{1}{(2-\beta)\lambda_{1}-\gamma}\big)\nonumber\\
&\quad+\big(\gamma^{(1-\alpha)/\alpha}+\frac{1}{\lambda_{1}p-\gamma}\big)^{p/2}\big]<1.
\end{align*}
Then there exists $\lambda_{L}>0$ such that for any $\lambda\in (0,\lambda_{L})$,$u_{\lambda}\in B_{p,\gamma}$.  The proof is complete. $\Box$ \\ \textbf{Remark 3.8.}  Fix some $x\in (0,L)$ and define the upper $p$th-moment Liapounov exponent $\overline{\gamma}(p)$ of $u$ as following (see \cite{FK})
\begin{align*}
\overline{\gamma}(p):=\lim_{t\rightarrow \infty}\sup\frac{1}{t}\ln\mathbb{E}|u(t,x)|^{p}, \ \mbox{for all} \ p\in(0,\infty).
\end{align*}
$u$ is called weakly intermittent, if $\overline{\gamma}(p)\in (0,\infty)$, for all $p\geq 2$. From theorem 3.7, it follows that the solution $u_{\lambda}$
is not of weak intermittent.

\section{\textbf{Noise excitation index of $p$th energy of the solutions to \emph{Eq.} (\ref{sfh})}}
In this section, we first prove that the second moment of the solution to \emph{Eq.} (\ref{sfh}) has a lower bound on a closed subinterval of $(0,L)$, and the second moment grows at most exponentially on $[0,L]$.
We then show that the $p$th moment of $u_{\lambda}(t,x)$ grows faster than a Mittag-Leffler function.
At last, we prove that the excitation index of the mild solution of  \emph{Eq.} (\ref{sfh}) is $\frac{2\alpha}{\alpha-1}$.\\
\textbf{Theorem 4.1.} 
Fix $\mu\in (0,L/2)$ and suppose $\sigma$ satisfies Assumption 1.2, then there exists constants $\kappa_{1},\kappa_{2}>0$ such that for all $t>0$,
\begin{align}\label{kappa}
\inf_{x\in [\mu,L-\mu]}\mathbb{E}|u_{\lambda}(t,x)|^{2}
\geq  \kappa_{1}E_{1-1/\alpha}(\lambda^{2}l^{2}_{\sigma}\kappa_{2}t^{(\alpha-1)/\alpha}),
\end{align}
where $E(\cdot)$ is the Mittag-Leffler functions (see formula (1.66) in \cite{Igor}),  defined by
\begin{align}\label{ML}
 E_{\beta}(z)=\sum_{n=0}^{\infty}\frac{z^{n}}{\Gamma(n\beta+1)}, \ \beta>0, \ z\in \mathbb{C},
\end{align}
where $\mathbb{C}$ denotes the complex plane. \\
\textbf{Proof.} From (i) of Lemma 3.5,  it follows that
\begin{align*}
p_{D}(t,x,y)\geq c_{0} \big(1\wedge \frac{\delta_{D}(x)^{\alpha/2}}{\sqrt{t}}\big)\big(1\wedge \frac{\delta_{D}(y)^{\alpha/2}}{\sqrt{t}}\big)\big(t^{-1/\alpha}\wedge \frac{t}{|x-y|^{1+\alpha}}\big).
\end{align*}
This together with (\ref{c1}) yield
\begin{align}\label{Foondun}
p_{D}(t,x,y)\geq c_{0}p(t,x,y), \ \mbox{for} \ 0<t\leq\mu^{\alpha}, \ x,y\in [\mu,L-\mu].
\end{align}
Since $|\sigma(u)|\geq  l_{\sigma}|u|$, we have
\begin{align}\label{lamdal}
&\mathbb{E}\big|\int_{0}^{t}\int_{0}^{L}p_{D}(t-s,x,y)\lambda\sigma(u_{\lambda}(s,y))w(dsdy)\big|^{2}\nonumber\\
&=\int_{0}^{t}\int_{0}^{L}p^{2}_{D}(t-s,x,y)\mathbb{E}|\lambda\sigma(u_{\lambda}(s,y))|^{2}dyds\nonumber\\
&\geq \lambda^{2}l^{2}_{\sigma}\int_{0}^{t}\int_{0}^{L}p^{2}_{D}(t-s,x,y)\mathbb{E}|u_{\lambda}(s,y)|^{2}dyds
\end{align}
Set
\begin{align}\label{at}
a(t)=\inf_{x\in[\mu,L-\mu]}\mathbb{E}|u_{\lambda}(t,x)|^{2}.
\end{align}
By (\ref{Foondun}), we get
\begin{align}\label{ulamdel}
\int_{0}^{t}\int_{0}^{L}p^{2}_{D}(t-s,x,y)\mathbb{E}|u_{\lambda}(s,y)|^{2}dyds\geq\int_{0}^{t}\int_{\mu}^{L-\mu}p^{2}_{D}(t-s,x,y)a(s)dyds.
\end{align}
For $0<t\leq \mu^{\alpha}$, by (\ref{Foondun}), we obtain
\begin{align}\label{mualpha}
\int_{0}^{t}\int_{\mu}^{L-\mu}p^{2}_{D}(t-s,x,y)a(s)dyds\geq c_{0}\int_{0}^{t}\int_{\mu}^{L-\mu}p^{2}(t-s,x,y)a(s)dyds.
\end{align}
Set $A:=[\mu,L-\mu]\cap \{y:|y-x|\leq (t-s)^{1/\alpha}\}$, since $0\leq t-s\leq t<\mu^{\alpha}$, we have $|A|\geq (t-s)^{1/\alpha}$, where $|A|$ denotes the volume of $A$. By (\ref{c1}), We have
\begin{align}\label{mu}
\int_{\mu}^{L-\mu}p^{2}(t-s,x,y)dy &\geq c^{2}_{1}(\alpha)\int_{\mu}^{L-\mu}\big((t-s)^{-2/\alpha}\wedge \frac{(t-s)^{2}}{|x-y|^{2+2\alpha}}\big)dy\nonumber\\
&\geq c^{2}_{1}(\alpha)\int_{A}\big((t-s)^{-2/\alpha}\wedge \frac{(t-s)^{2}}{|x-y|^{2+2\alpha}}\big)dy\nonumber\\
&\geq c^{2}_{1}(\alpha)\int_{A}(t-s)^{-2/\alpha}dy\nonumber\\
&\geq c^{2}_{1}(\alpha)(t-s)^{-1/\alpha}.
\end{align}
By (\ref{ulamdel}), (\ref{mualpha}) and (\ref{mu}), for all $0<t\leq \mu^{\alpha}$, we have
\begin{align}\label{txye}
\int_{0}^{t}\int_{0}^{L}p^{2}_{D}(t-s,x,y)\mathbb{E}|u_{\lambda}(s,y)|^{2}dyds\geq c_{0}c^{2}_{1}(\alpha)\int_{0}^{t}(t-s)^{-1/\alpha}a(s)ds.
\end{align}
Set
\begin{align*}
(\mathcal{G}_{D}u_{\lambda})(t,x)=\int_{0}^{L}p_{D}(t,x,y)u_{0}(y)dy.
\end{align*}
For $0<t\leq \mu^{\alpha}$ and $x\in [\mu,L-\mu]$, we have
\begin{align*}
|(\mathcal{G}_{D}u_{\lambda})(t,x)|^{2}&=|\int_{0}^{L}p_{D}(t,x,y)u_{0}(y)dy|^{2}\nonumber\\
&\geq (\inf_{x\in [\mu,L-\mu]}u_{0}(x))^{2}(\int_{\mu}^{L-\mu}p_{D}(t,x,y)dy)^{2}\nonumber\\
&\geq c^{2}_{0}(\inf_{x\in [\mu,L-\mu]}u_{0}(x))^{2}\int_{\mu}^{L-\mu}p(t,x,y)dy\nonumber\\
&\geq c^{2}_{0}c_{2}(\int_{\mu}^{L-\mu}p(t,x,y)dy)^{2},
\end{align*}
where $c_{2}=(\inf_{x\in [\mu,L-\mu]}u_{0}(x))^{2}$. By similar argument to the proof of  (\ref{mu}), we have
\begin{align}\label{cons}
(\mathcal{G}_{D}u_{\lambda})(t,x)\geq c_{3},
\end{align}
where $c_{3}>0$ is a constant. \\
Taking the second moment to (\ref{Walsh}), we have
\begin{align*}
\mathbb{E}|u_{\lambda}(t,x)|^{2}=|(\mathcal{G_{D}}u_{\lambda})(x,t)|^{2}+\mathbb{E}\big|\int_{0}^{t}\int_{0}^{L}p_{D}(t-s,x,y)\lambda\sigma(u_{\lambda}(s,y))w(dsdy)\big|^{2}.
\end{align*}
By (\ref{lamdal}), (\ref{txye}) and (\ref{cons}), for $t\in (0,\mu^{\alpha}]$, we get
\begin{align}\label{integral}
\mathbb{E}|u_{\lambda}(t,x)|^{2}\geq c^{2}_{3}+\lambda^{2}l^{2}_{\sigma}c_{0}c^{2}_{1}(\alpha)\int_{0}^{t}(t-s)^{-1/\alpha}a(s)ds.
\end{align}
For any fixed $t,T>0$,  we get
\begin{align}\label{change}
&\mathbb{E}|u_{\lambda}(T+t,x)|^{2}\nonumber\\
&=|(\mathcal{G}_{D}u_{\lambda})(T+t,x)|^{2}+\lambda^{2}\int_{0}^{T+t}\int_{0}^{L}p^{2}_{D}(T+t-s)\mathbb{E}|\sigma(u_{\lambda}(s,y))|^{2}dyds\nonumber\\
&=|(\mathcal{G}_{D}u_{\lambda})(T+t,x)|^{2}+\lambda^{2}\int_{0}^{T}\int_{0}^{L}p^{2}_{D}(T+t-s)\mathbb{E}|\sigma(u_{\lambda}(s,y))|^{2}dyds\nonumber\\
&\quad+\lambda^{2}\int_{0}^{t}\int_{0}^{L}p^{2}_{D}(t-s)\mathbb{E}|\sigma(u_{\lambda}(T+s,y))|^{2}dyds.
\end{align}
Then
\begin{align}\label{similar}
&\mathbb{E}|u_{\lambda}(T+t,x)|^{2}\nonumber\\
&\geq |(\mathcal{G}_{D}u_{\lambda})(T+t,x)|^{2}|+\lambda^{2}l^{2}_{\sigma}\int_{0}^{t}\int_{0}^{L}p^{2}_{D}(t-s)\mathbb{E}|u_{\lambda}(T+s,y))|^{2}dyds.
\end{align}
Then we see that for any $t>0$, (\ref{integral}) holds.
From the definition of $a(t)$, it follows that
\begin{align}\label{equality}
a(t)\geq c^{2}_{3}+\lambda^{2}l^{2}_{\sigma}c_{4}\int_{0}^{t}(t-s)^{-1/\alpha}a(s)ds, \ t>0,
\end{align}
where $c_{4}=c_{0}c^{2}_{1}(\alpha)$. \\
By Lemma 2.4, we have
\begin{align}\label{at}
a(t)\geq c^{2}_{3}F_{1-1/\alpha}((\lambda^{2}l^{2}_{\sigma}c_{4}\Gamma(1-1/\alpha))^{\alpha/(\alpha-1)}t).
\end{align}
From Lemma 2.4 and (\ref{ML}), We see that
$F_{\beta}(z)=E_{\beta}(z^{\beta})$ for any $\beta>0$. Then by (\ref{at}), we have
\begin{align}\label{ef}
a(t)&\geq c^{2}_{3}E_{1-1/\alpha}[((\lambda^{2}l^{2}_{\sigma}c_{4}\Gamma(1-1/\alpha))^{\alpha/(\alpha-1)}t)^{(\alpha-1)/\alpha}]\nonumber\\
&=c^{2}_{3}E_{1-1/\alpha}(\lambda^{2}l^{2}_{\sigma}c_{4}\Gamma(1-1/\alpha)t^{(\alpha-1)/\alpha}).
\end{align}
The proof is complete. $\Box$\\
\textbf{Theorem 4.2.} There exists constants $\kappa_{3}>0$, $\kappa_{4}>0$ such that for all $t>0$,
\begin{align*}
\sup_{x\in [0,L]}\mathbb{E}|u_{\lambda}(t,x)|^{2}\leq \kappa_{3}\exp\big(\kappa_{4}(\lambda^{2}L^{2}_{\sigma})^{\alpha/(\alpha-1)}t\big),
\end{align*}
where $u_{\lambda}(t,x)$ is the unique mild solution of (\ref{sfh}). \\
\textbf{Proof.}  Set
\begin{align*}
U(t)=\sup_{x\in [0,L]}\mathbb{E}|u_{\lambda}(t,x)|^{2},
\end{align*}
and \begin{align*}
(\mathcal{G}_{D}u_{\lambda})(t,x)=\int_{0}^{L}p_{D}(t,x,y)u_{0}(y)dy.
\end{align*}
We have
\begin{align}\label{Walsh}
\mathbb{E}|u_{\lambda}(t,x)|^{2}=|(\mathcal{G_{D}}u_{\lambda})(x,t)|^{2}+\mathbb{E}\big|\int_{0}^{t}\int_{0}^{L}p_{D}(t-s,x,y)\lambda\sigma(u_{\lambda}(s,y))w(dsdy)\big|^{2}.
\end{align}
By (\ref{inequality}),
\begin{align}\label{mathcal}
|(\mathcal{G}_{D}u_{\lambda})(t,x)|^{2}&\leq|\int_{0}^{L}p_{D}(t,x,y)u_{0}(y)dy|^{2}\nonumber\\
&\leq c_{1}|\int_{0}^{L}p_{D}(t,x,y)dy|^{2}\nonumber\\
&\leq c_{1}|\int_{0}^{L}p(t,x,y)dy|^{2}\nonumber\\
&\leq c_{1}.
\end{align}
By (\ref{inequality}) and (\ref{c1}), we have
\begin{align}\label{math}
&\mathbb{E}\big|\int_{0}^{t}\int_{0}^{L}p_{D}(t-s,x,y)\lambda\sigma(u_{\lambda}(s,y))w(dsdy)\big|^{2}\nonumber\\
&\leq \lambda^{2}L^{2}_{\sigma}\int_{0}^{t}\int_{0}^{L}p^{2}_{D}(t-s,x,y)\mathbb{E}|u_{\lambda}(s,y)|^{2}dyds\nonumber\\
&\leq \lambda^{2}L^{2}_{\sigma}\int_{0}^{t}U(s)\big(\int_{0}^{L}p^{2}_{D}(t-s,x,y)dy\big)ds\nonumber\\
&\leq \lambda^{2}L^{2}_{\sigma}\int_{0}^{t}U(s)p_{D}(2(t-s),x,x)ds\nonumber\\
&\leq c_{2}\lambda^{2}L^{2}_{\sigma}\int_{0}^{t}(t-s)^{-1/\alpha}U(s)ds.
\end{align}
From (\ref{Walsh}), (\ref{mathcal}) and (\ref{math}), it follows that
\begin{align*}
U(t)\leq c_{1}+c_{2}\lambda^{2}L^{2}_{\sigma}\int_{0}^{t}(t-s)^{-1/\alpha}U(s)ds.
\end{align*}
By Lemma 7.1.1 in \cite{DH}, we have
\begin{align}\label{U}
U(t)\leq c_{1}F_{1-1/\alpha}((c_{2}\lambda^{2}L^{2}_{\sigma}\Gamma(1-1/\alpha))^{\alpha/(\alpha-1)}t),
\end{align}
where $F$ is the same function as that in Lemma 2.4.

It is known that (see \cite{Bazhlekova}, formula (2.9)) the Mittag-Leffler function $E_{\tau}(\omega t^{\tau})$ statistics the following inequality
\begin{align}\label{tau}
E_\tau(\omega t^{\tau})\leq c\exp(\omega ^{1/\tau}t), \ t\geq 0, \ \tau\in(0,2).
\end{align}
Note that $F_{\beta}(z)=E_{\beta}(z^{\beta})$ for any $\beta>0$. By (\ref{U}) and (\ref{tau}), we obtain the result.
\ $\Box$ \ \\
\textbf{Theorem 4.3.} Under the assumptions in Theorem 4.1,  there exists $\lambda_{0}> \lambda_{L}$ sucht that for all $\lambda\in(\lambda_{0},\infty)$
and $x\in [\mu,l-\mu]$,
\begin{align}\label{energy}
(\kappa_{2}\lambda^{2}l^{2}_{\sigma})^{\alpha/(\alpha-1)}<\lim_{t\rightarrow \infty}\inf\frac{1}{t}\log\mathbb{E}|u_{\lambda}(t,x)|^{p}\leq \lim_{t\rightarrow \infty}\sup\frac{1}{t}\log\mathbb{E}|u_{\lambda}(t,x)|^{p}_{L^{\infty}}<\infty,
\end{align}
where $\kappa_{2}$ is the same constant as that in Theorem 4.1. \\
\textbf{Proof.} By Theorem 4.1, for all $t>0$,
\begin{align*}
\inf_{x\in [\mu,L-\mu]}\mathbb{E}|u_{\lambda}(t,x)|^{2})\geq \kappa_{1}E_{1-1/\alpha}(\lambda^{2}l^{2}_{\sigma}\kappa_{2}t^{(\alpha-1)/\alpha}).
\end{align*}
The Mittag-Leffler functions  has the following asymptotic property as $z\rightarrow \infty$ (see Theorem 1.3 in \cite{Igor}): if $0<\tau<2$, $\mu$ is an arbitrary real number such that $\pi \tau/2<\mu<\min\{\pi,\pi\tau\}$, then for arbitrary integer $q\geq 1$, the following expansion holds:
\begin{align}\label{hold}
E_{\tau}(z)=\frac{1}{\tau}\exp(z^{1/\tau})-\sum_{k=1}^{q}\frac{z^{-k}}{\Gamma(1-\tau k)}+O(|z|^{-1-q}),
\end{align}
$|z|\rightarrow \infty$, $|\arg(z)|\leq \mu$, where $\arg(z)$ denotes the principal value of the argument of $z$.
Choosing $q=1$, by (\ref{hold}), we get
\begin{align}\label{arg}
E_{\tau}(z)=\frac{1}{\tau}\exp(z^{1/\tau})-\frac{z^{-1}}{\Gamma(1-\tau)}+O(|z|^{-2}),
\end{align}
$|z|\rightarrow \infty$, $|\arg(z)|\leq \mu$. \\
By (\ref{kappa}), we have
\begin{align}\label{compare}
&\inf_{x\in[\mu,L-\mu]}\mathbb{E}|u_{\lambda}(t,x)|^{2}\nonumber\\
&\geq \frac{\alpha}{\alpha-1}\kappa_{1}\exp((\lambda^{2}l^{2}_{\sigma}\kappa_{2})^{\alpha/(\alpha-1)}t)-\frac{t^{(1-\alpha)/\alpha}}{\Gamma(1/\alpha)\lambda^{2}l^{2}_{\sigma}\kappa_{2}}\nonumber\\
&\quad+O(\lambda^{-4}), \ \mbox{as} \ \lambda\rightarrow \infty.
\end{align}
We can choose $\lambda_{0}>\lambda_{L}$ sufficiently large such that for all $\lambda>\lambda_{0}$,
\begin{align}\label{sufficient}
\inf_{x\in[\mu,L-\mu]}\mathbb{E}|u_{\lambda}(t,x)|^{2}\geq c_{5}\exp((\lambda^{2}l^{2}_{\sigma}\kappa_{2})^{\alpha/(\alpha-1)}t),
\end{align}
where $c_{5}>0$ ia a constant.
Then
\begin{align}\label{this}
\lim_{t\rightarrow \infty}\inf\frac{1}{t}\log(\inf_{x\in [\mu,L-\mu]}\mathbb{E}|u_{\lambda}(t,x)|^{2})\geq (\kappa_{2}\lambda^{2}l^{2}_{\sigma})^{\alpha/(\alpha-1)}.
\end{align}
By Jensen's inequality, for $p\geq2$, we have
\begin{align*}
(\mathbb{E}|u_{\lambda}(t,x)|^{2})^{1/2}\leq (\mathbb{E}|u_{\lambda}(t,x)|^{p})^{1/p}.
\end{align*}
This together with  Theorem 3.4 yield that the inequality (\ref{energy}) holds.  \ $\Box$ \ \\
\textbf{Remark 4.4.} By the definition of weakly intermittency, we see that the solution of \emph{Eq.}(\ref{sfh}) is weak intermittent under the assumptions in Theorem 3.7.

Under the assumptions in Theorem 3.7, for $\lambda\in (0,\lambda_{L})$, the $p$th energy of the solution of \emph{Eq.}(\ref{sfh}) satisfies the following inequality
\begin{align}\label{energy1}
-\infty<\lim_{t\rightarrow \infty}\sup\frac{1}{t}\log\Phi_{p}(t,\lambda)<0.
\end{align}
Since for $p\geq 2$, $\mathbb{E}\|u_{\lambda}(t)\|^{p}_{L^{p}}\leq \mathbb{E}\|u_{\lambda}(t)\|^{p}_{L^{\infty}}$, by Theorem 3.7,
it is easy to see that (\ref{energy1}) holds.  \\
\textbf{Corollary 4.5.} Under the assumptions in Theorem 4.1 and Theorem 3.4, then for $\lambda>\lambda_{0}$, the $p$th energy of the solution of \emph{Eq.}(\ref{sfh}) satisfies
\begin{align}\label{energy2}
\frac{1}{2}(\kappa_{2}\lambda^{2}l^{2}_{\sigma})^{\alpha/(\alpha-1)}<\lim_{t\rightarrow \infty}\inf\frac{1}{t}\log\Phi_{p}(t,\lambda)<\infty,
\end{align}
where $\kappa_{2}$ is the same constant as that in Theorem 4.1. \\
\textbf{Proof.} For $p\geq2$ and $\lambda>\lambda_{0}$, by Jensen's inequality and (\ref{sufficient}),
\begin{align*}
\mathbb{E}\int_{0}^{L}|u_{\lambda}(t,x)|^{p}dx&\geq \int_{0}^{1}(\mathbb{E}|u_{\lambda}(t,x)|^{2})^{p/2}dx\\
&\geq \int_{\mu}^{L-\mu}(\inf_{x\in [\mu,L-\mu]}\mathbb{E}|u_{\lambda}(t,x)|^{2})^{p/2}\nonumber\\
&\geq (L-2\mu)(\inf_{x\in [\mu,L-\mu]}\mathbb{E}|u_{\lambda}(t,x)|^{2})^{p/2}\nonumber\\
&\geq (L-2\mu)\big(c_{5}\exp((\lambda^{2}l^{2}_{\sigma}\kappa_{2})^{\alpha/(\alpha-1)}t))^{p/2}.
\end{align*}
Then
\begin{align}\label{corollary}
\log\Phi_{p}(t,\lambda)&=\frac{1}{p}\log(\mathbb{E}\|u_{\lambda}(t)\|^{p}_{L^{p}})\nonumber\\
&=\frac{1}{p}\log\big(\mathbb{E}\int_{0}^{L}|u_{\lambda}(t,x)|^{p}dx\big)\nonumber\\
&\geq \frac{1}{p}\log(L-2\mu)+\frac{1}{2}\log\big(c_{5}\exp((\lambda^{2}l^{2}_{\sigma}\kappa_{2})^{\alpha/(\alpha-1)}t)\big).
\end{align}
Therefore, we get
\begin{align}
\lim_{t\rightarrow \infty}\inf\frac{1}{t}\log\Phi_{p}(t,\lambda)\geq \frac{1}{2}(\kappa_{2}\lambda^{2}l^{2}_{\sigma})^{\alpha/(\alpha-1)}.
\end{align}
This together with Theorem 3.4 yield that (\ref{energy2}).  \ $\Box$ \\
\textbf{Lemma 4.6.} (Minkowski's Inequality) Let $X$ and $Y$ be random variables. Then, for $1\leq p<\infty$,
\begin{align*}
(\mathbb{E}|X+Y|^{p})^{1/p}\leq (\mathbb{E}|X|^{p})^{1/p}+(\mathbb{E}|Y|^{p})^{1/p}.
\end{align*}
\textbf{Theorem 4.7.} Under the assumptions in Theorem 4.1, there exists a constant $c_{p}>0$ such that for all $t>0$,
\begin{align}\label{assumption}
c_{p}l^{2\alpha/(\alpha-1)}_{\sigma}t&\leq \lim_{\lambda\rightarrow \infty}\inf\lambda^{2\alpha/(1-\alpha)}\log\big(\inf_{x\in [\mu,L-\mu]}\mathbb{E}|u_{\lambda}(t,x)|^{p}\big)\nonumber\\
&\leq \lim_{\lambda\rightarrow \infty}\sup\lambda^{2\alpha/(1-\alpha)}\log\big(\sup_{x\in [0,L]}\mathbb{E}|u_{\lambda}(t,x)|^{p}\big)\leq c^{-1}_{p}L^{2\alpha/(\alpha-1)}_{\sigma}t.
\end{align}
\textbf{Proof.} First, we prove the upper bound. For all $t>0$, since $u_{0}$ is continuous on $[0,L]$,
\begin{align*}
\sup_{x\in [0,L]}|(\mathcal{G}_{D}u_{\lambda})(t,x)|=\sup_{x\in [0,L]}\int_{0}^{L}|p_{D}(t,x,y)u_{0}(y)|dy\leq c_{1}.
\end{align*}
By similar arguments as that in Lemma 3.2, we have
\begin{align}\label{Henry}
&\big(\mathbb{E}\big|\int_{0}^{t}\int_{0}^{L}p_{D}(t-s,x,y)\lambda\sigma(u_{\lambda}(s,y))w(dsdy)\big|^{p}\big)^{2/p}\nonumber\\
&\leq c_{2}\lambda^{2}L^{2}_{\sigma}\int_{0}^{L}\int_{0}^{t}p^{2}_{D}(t-s,x,y)(\mathbb{E}|u_{\lambda}(s,y)|^{p})^{2/p}dyds\nonumber\\
&\leq c_{2}\lambda^{2}L^{2}_{\sigma}\int_{0}^{L}\int_{0}^{t}p^{2}_{D}(t-s,x,y)(\sup_{y\in [0,L]}\mathbb{E}|u_{\lambda}(s,y)|^{p})^{2/p}dyds\nonumber\\
&\leq c_{2}\lambda^{2}L^{2}_{\sigma}L\int_{0}^{t}p^{2}_{D}(t-s,x,y)(\sup_{y\in [0,L]}\mathbb{E}|u_{\lambda}(s,y)|^{p})^{2/p}ds\nonumber\\
&\leq c_{2}\lambda^{2}L^{2}_{\sigma}L\int_{0}^{t}p_{D}(2(t-s),x,x)(\sup_{y\in [0,L]}\mathbb{E}|u_{\lambda}(s,y)|^{p})^{2/p}ds\nonumber\\
&\leq c_{3}\lambda^{2}L^{2}_{\sigma}L\int_{0}^{t}(t-s)^{-1/\alpha}(\sup_{y\in [0,L]}\mathbb{E}|u_{\lambda}(s,y)|^{p})^{2/p}.
\end{align}
By Lemma 4.6,
\begin{align}\label{estimate}
(\sup_{x\in [0,L]}\mathbb{E}|u_{\lambda}(t,x)|^{p})^{2/p}&\leq 2(\mathbb{E}\sup_{x\in [0,L]}|(\mathcal{G}_{D}u_{\lambda})(t,x)|^{p})^{2/p}\nonumber\\
&\quad+2(\mathbb{E}\big|\int_{0}^{t}\int_{0}^{L}p_{D}(t-s,x,y)\lambda\sigma(u_{\lambda}(s,y))w(dsdy)\big|^{p})^{2/p}\nonumber\\
&\leq c_{4}+c_{5}\lambda^{2}L^{2}_{\sigma}L\int_{0}^{t}(t-s)^{-1/\alpha}(\sup_{x\in [0,L]}\mathbb{E}|u_{\lambda}(s,x)|^{p})^{2/p}ds.
\end{align}
where $c>0$ is a constant. \\
By Lemma 7.1.1 in \cite{DH} and (\ref{estimate}), for any $t>0$,
\begin{align}\label{continuous}
(\sup_{x\in [0,L]}\mathbb{E}|u_{\lambda}(t,x)|^{p})^{2/p}&\leq c_{4}F_{1-1/\alpha}((c_{5}\lambda^{2}L^{2}_{\sigma}L\Gamma(1-1/\alpha))^{\alpha/(\alpha-1)}t)\nonumber\\
&=c_{4}E_{1-1/\alpha}(c_{5}\lambda^{2}L^{2}_{\sigma}L\Gamma(1-1/\alpha)t^{(\alpha-1)/\alpha})\nonumber\\
&\leq c_{6}\exp[(c_{5}\lambda^{2}L^{2}_{\sigma}L\Gamma(1-1/\alpha))^{\alpha/(\alpha-1)}t].
\end{align}
Therefore, we obtain the upper bound in (\ref{assumption}). \\
By Theorem 4.1, for any $t>0$,
\begin{align}
\log(\inf_{x\in [\mu,L-\mu]}\mathbb{E}|u_{\lambda}(t,x)|^{2})\geq \log (\kappa_{1}E_{1-1/\alpha}(\lambda^{2}l^{2}_{\sigma}\kappa_{2}t^{(\alpha-1)/\alpha})).
\end{align}
Then from (\ref{sufficient}), for $\lambda>\lambda_{0}$, it follows that
\begin{align}\label{Jensen}
\log(\inf_{x\in [\mu,L-\mu]}\mathbb{E}|u_{\lambda}(t,x)|^{2})&\geq c_{5}\log(\exp((\lambda^{2}l^{2}_{\sigma}\kappa_{2})^{\alpha/(\alpha-1)}t))\nonumber\\
&=c_{5}(\lambda^{2}l^{2}_{\sigma}\kappa_{2})^{\alpha/(\alpha-1)}t.
\end{align}
By Jensen's inequality and (\ref{Jensen}), we obtain the lower bound in  (\ref{assumption}). \ $\Box$  \

Set 
\begin{align*}
f_{p,t}(\lambda):=\left(\sup_{x\in [0,L]}\mathbb{E}|u_{\lambda}(t,x)|^{p}\right)^{1/p}.
\end{align*}
\textbf{Lemma 4.8} Fix $t>0$, then
\begin{align}\label{loglog}
\lim_{\lambda\rightarrow \infty}\frac{\log\log f_{p,t}(\lambda)}{\log \lambda} \leq\frac{2\alpha}{\alpha-1}.
\end{align}
\textbf{Proof.}
By (\ref{continuous}),
\begin{align*}
(\sup_{x\in [0,L]}\mathbb{E}|u_{\lambda}(t,x)|^{p})^{1/p}\leq \big(c_{6}\exp[(c_{5}\lambda^{2}L^{2}_{\sigma}L\Gamma(1-1/\alpha))^{\alpha/(\alpha-1)}t]\big)^{1/2}.
\end{align*}
It is easy to see that (\ref{loglog}) holds. \ $\Box$

For any $\mu\in (0,L/2)$, set
\begin{align*}
I_{p,\mu,t}(\lambda)=\left(\inf_{x\in [\mu,L-\mu]}\mathbb{E}|u_{\lambda}(t,x)|^{p}\right)^{1/p}.
\end{align*}
\textbf{Lemma 4.9.} For any $\mu\in (0,L/2)$, fix $t>0$,
\begin{align*}
\lim_{\lambda\rightarrow \infty}\frac{\log\log I_{p,\mu,t}(\lambda)}{\log \lambda} \geq\frac{2\alpha}{\alpha-1}.
\end{align*}
\textbf{Proof.}
By (\ref{sufficient}), for all $\lambda>\lambda_{0}$, we have
\begin{align}\label{acf}
\inf_{x\in[\mu,L-\mu]}\mathbb{E}|u_{\lambda}(t,x)|^{2}\geq c_{5}\exp((\lambda^{2}l^{2}_{\sigma}\kappa_{2})^{\alpha/(\alpha-1)}t)
\end{align}
By Jensen's inequality, for $p\geq2$,
\begin{align}\label{Jensen Inequality}
(\mathbb{E}|u_{\lambda}(t,x)|^{2})^{1/2}\leq (\mathbb{E}|u_{\lambda}(t,x)|^{p})^{1/p}.
\end{align}
By (\ref{acf}) and (\ref{Jensen Inequality}), we obtain the result.  \ $\Box$  \ \\
\textbf{Theorem 4.10.} The excitation index of the solution of (\ref{sfh}) is $\frac{2\alpha}{\alpha-1}$. \\
\textbf{Proof.} Since the $p$th energy is defined by
\begin{align*}
\Phi_{p}(t,\lambda)=(\mathbb{E}\|u_{\lambda}(t)\|^{p}_{L^{p}})^{1/p}=\big(\int_{0}^{L}\mathbb{E}|u_{\lambda}(t,x)|^{p}dx\big)^{1/p}, \ t>0.
\end{align*}
By Lemma 4.8, Lemma 4.9 and the following inequalities
\begin{align*}
\int_{0}^{L}\mathbb{E}|u_{\lambda}(t,x)|^{p}dx\leq L\sup_{x\in [0,L]}\mathbb{E}|u_{\lambda}(t,x)|^{p},
\end{align*}
\begin{align*}
\int_{0}^{L}\mathbb{E}|u_{\lambda}(t,x)|^{p}dx\geq (L-2\mu)\inf_{x\in [\mu,L-\mu]}\mathbb{E}|u_{\lambda}(t,x)|^{p},
\end{align*}
we obtain the result. \ $\Box$  \ \\
\textbf{Remark 4.11.} For stochastic heat equations, when $p=2$, Foondun and Joseph \cite{FJ} showed that $e_{2}=4$, when $p>2$, Xie \cite{Xiebin}
proved that $e_{p}=4$. If we consider the case $\alpha\rightarrow 2$, the stochastic fractional heat equations reduces to stochastic heat equations, then we obtain $e_{p}\rightarrow 4$.

\section{Acknowledgements}
This author is supported by National Natural Science Foundation of China under the contract
No.11571269, China Postdoctoral Science Foundation Funded Project under contract Nos.2015M572539, 2016T90899 and Shaanxi Province Natural Science Research Project under the contract 2017JM1011.



\end{document}